# Vortices and atoms in the Maxwellian era

Isobel Falconer

School of Mathematics and Statistics,
University of St Andrews,
North Haugh,
St Andrews,
KY16 9SS

Email: ijf3@st-andrews.ac.uk

## Abstract
The mathematical study of vortices began with Herman von Helmholtz's pioneering study in 1858. It was pursued vigorously over the next two decades, largely by British physicists and mathematicians, in two contexts: Maxwell's vortex analogy for the electromagnetic field and William Thomson's (Lord Kelvin) theory that atoms were vortex rings in an all-pervading ether. By the time of Maxwell's death in 1879, the basic laws of vortices in a perfect fluid in three-dimensional Euclidean space had been established, as had their importance to physics. Early vortex studies were embedded in a web of issues spanning the fields we now know as "mathematics" and "physics" – fields which had not yet become institutionally distinct disciplines but overlapped. This paper investigates the conceptual issues with ideas of force, matter, and space, that underlay mechanics and led to vortex models being an attractive proposition for British physicists, and how these issues played out in the mathematics of vortices, paying particular attention to problems around continuity. It concludes that while they made valuable contributions to hydrodynamics and the nascent field of topology, the British ultimately failed in their more physical objectives.

## Keywords
vortices, vortex atoms, Maxwell, William Thomson, Kelvin, Helmholtz

## 1. Introduction
The mathematical study of vortices began with Herman von Helmholtz's pioneering study in 1858, translated by Peter Guthrie Tait in 1867 [1,2]. It was pursued vigorously over the next two decades, largely by British physicists and mathematicians with whom Helmholtz was on close terms. This paper examines the work of an influential early group, whom Smith [3] has dubbed the "North British".[1] The group was led by William Thomson (referred to throughout the rest of this paper by his later title, Kelvin) and his friends Peter Guthrie Tait and James Clerk Maxwell. All three were raised in Scotland, concluded their education by studying mathematics at Cambridge, and retained a Scottish base throughout their lives.[2] By the time of Maxwell's premature death in 1879, understanding of vortices was, "without doubt, the most important addition to the theory of fluid motion" of the past few decades [7 p63]. It was an important field of enquiry, taken up enthusiastically by a slightly younger, Cambridge-trained, group of mathematical physicists, notably William Hicks, Horace Lamb, and Joseph John Thomson (referred to as Thomson). It was the main topic of Hicks' "Report on Recent Progress in Hydrodynamics" commissioned by the British Association for the Advancement of Science [7], and the subject of the Adams Prize set by Cambridge for 1882, "A general investigation of the action upon each other of two closed vortices in a perfect incompressible fluid," won by Thomson [8].

---

[1] For works encompassing German and French contributions more fully, see Darrigol [4], Meleshko, Gourjii and Krasnopolskaya [5], and Frisch, Grimberg and Villone [6]

To understand the way the subject developed in these early years, it is important to investigate the reciprocal interactions between mathematical conceptions and methods, physical or empirical considerations, and social context. The remainder of this Introduction gives a brief overview of the context. The following section looks at the issues with mechanics that led to vortex models being an attractive proposition for the North British. The subsequent sections examine how these issues played out in the mathematics of vortices, paying particular attention to problems around continuity.

Early vortex studies were embedded in a web of issues spanning the fields we now know as "mathematics" and "physics" – fields which had not yet become institutionally distinct disciplines but overlapped, with shifting boundaries. On the more mathematical side the nature of the symbols used in algebra, and approaches to the foundations of the calculus, were debated, with different attitudes generally observable on either side of the English Channel. From the more physical perspective, the nature of space, time, and mass, the status of Newton's laws, and the role of such concepts as force and energy, were all disputed. Challenging the union of the two was the fundamental issue of the relationship between mathematics and the physical world as differing assumptions were made about the interface between empirically observable and unobservable entities, and as the development of non-Euclidean geometries exposed the sensory basis of mathematical truth claims.

In this context, the hydrodynamics – particularly the vortex dynamics - of an all-pervading but insensible fluid, the ether, appeared to offer a promising way forward to the North British physicists. Unlike their Continental counterparts, their interest was neither in the development of pure mathematical methods, nor in the increasingly powerful mathematical description of observed flows described by Darrigol [4]. Instead, as Duhem remarked, "Whereas the French or German physicist intends the algebraic part of a theory to replace just the series of syllogisms used to develop this theory, the English physicist regards the algebra as playing the part of a *model*." [9 p79]. The models were analogues of physical systems. They might, but more usually did not, describe some real underlying physical process, but at the level of macroscopic predictions they matched observation. Thus, a single mathematical theory might be used to model very different physical phenomena, as when Kelvin used fluid flow to model both heat and electromagnetism [10]. In such cases, understanding of one phenomenon might inform understanding of the other, with the mathematical model as an intermediary.

---

[2] Smith argues that the "North British" group were "not assimilable to Cambridge Anglicanism" [3 p7], and hence did not include George Gabriel Stokes, whose interests in hydrodynamics were rather different as discussed in section 4.

Ether theories had been around for a long time but rose to prominence in the 19[th] century [11]. They coalesced around the resurgence of the wave theory of light following Thomas Young's and Augustin Fresnel's work on diffraction, Michael Faraday's on electromagnetic induction, and Bernhard Reimann's attempts to unify light, gravitation, and heat in a single mathematical ether theory founded on differential equations. Two, closely connected, strands are evident in the British pursuit of ether vortex models. The first was electromagnetism; the second a theory of matter. In the early 1860s, inspired by Kelvin's suggestion that magnetism must have some sort of rotational origin, Maxwell famously extended Kelvin's fluid flow analogy for electrostatics and used a "sea" of ether vortices to model Faraday's lines of force and derive the equations of electromagnetism [12]. Although Maxwell subsequently discarded this heuristic model, he retained the mathematical analogy between the equations of the electromagnetic field and vortex dynamics [13–15]. In 1867, encouraged by Maxwell's success with vortex models, Kelvin exploited Helmholtz's observation that vortex rings in an incompressible frictionless fluid would be permanent, and could not be created or destroyed [1]. Kelvin suggested that molecules or atoms might be vortex rings in the ether. Atoms were unobservable, and many viewed them as purely hypothetical, but they had become increasingly important as an explanation for both the macroscopic properties of gases and the spectra of the elements. Ethereal vortices provided a way in which microscopic discrete entities which maintained a stable identity – atoms – could be viewed as epiphenomena of the continuous medium, and hence analysed using calculus-based methods.

Kelvin's vortex atom was enthusiastically developed by over a dozen British mathematicians and physicists, while Maxwell's death in 1879 left a number of unresolved problems and ambiguities in his electromagnetic theory that his followers struggled to understand. The resultant physical models have been extensively discussed, especially by [14,16–18]. To understand the enthusiasm, we investigate, first, debates over the foundations of mechanics and the relationship between mathematics and the physical world that had been festering for the previous 200 years.

## 2. Debates over foundations of mechanics

Development of vortex theory in the second half of the nineteenth century took place as mathematicians and physicists attempted to resolve fundamental issues around concepts of "force" and "space", and a disjoint between the mathematics of a continuum and the finite extension of matter. These problems, that shaped the mathematics of vortices, are outlined in this section.

Dating back to the seventeenth century, and the work of philosophers such as Descartes and Newton, most nineteenth century men of science viewed classical mechanics as the fundamental basis of physics: the ultimate goal was a mechanical explanation of natural phenomena. Kelvin stated an extreme version of this aim in his Baltimore Lectures of 1884: "It seems to me that the test of 'Do we or not understand a particular subject in physics?' is,

'Can we make a mechanical model of it?'" [19 vol2 p830]. Although Kelvin's Continental colleagues did not share the British pre-occupation with modelling as the means to a mechanical explanation, they shared the ultimate goal. Helmholtz, for example, aimed to "refer natural phenomena back to unchangeable attractive and repulsive forces, whose intensity depends solely upon distance," [20 p117, 21,22]. His protégé, Heinrich Hertz, concurred: "All physicists agree that the problem of physics consists in tracing the phenomena of nature back to the simple laws of mechanics," [23 p38].

Generally, such mechanical reductionism entailed explaining non-mechanical phenomena, such as light or heat, by the motion of unobservable – usually microscopic – mechanical systems. However, two philosophical problems challenged these attempts from the start. First was the question of whether the validity of the principles of mechanics rested on empirical experience (Newton), on rational reflection (Descartes), or on a Kantian combination of the two whereby a-priori concepts shape our empirical experience. Second was the issue of the applicability to the real world of a mathematical mechanics that had developed through the eighteenth century on the basis of an idealized world. Both issues were exemplified in debates over ideas of force and action at a distance.

Ideas about force in mechanics derived from Newton's second law. However, the status of the law – and of Newton's other two laws – was unclear. Was the second law a law of nature derived from the real world and positing force as the cause of motion, or was it a definition of force? In the eighteenth century D'Alembert, for example, took the latter view; for him *force* was defined as *mass* x *acceleration.* While this might be satisfactory mathematically, it led to the logical physical question: how could a mathematical definition act as a physical law and determine the motion of bodies? However, the former view, of force as the cause of a change of motion, was hardly more satisfactory. It was held by Euler but was rejected by those such as d'Alembert as being overly metaphysical [23]. In the nineteenth century such metaphysical arguments led Helmholtz into difficulties. He initially attempted to identify force as the final cause of the motion of matter that constituted physical explanations [10,20,21]. Being a fundamental cause, forces must be an invariable function of the distances between the material points on which they act, and hence Helmholtz argued that they must be central. His further "proof" that non-central forces, such as the velocity- and acceleration-dependent electric force suggested by Wilhelm Weber, would permit perpetual motion was demolished by Rudolph Lipschitz and Rudolph Clausius [21]. The episode prompted Helmholtz to withdraw from his metaphysical stance and adopt a more empiricist position, viewing the aim of science as the discovery of the laws relating appearances, and force as merely an assertion of law-like behaviour [24]. He was not alone in moving towards empiricism, a move that was prompted also by growing awareness of Riemann's claims that Euclidean geometry had only empirical, rather than universal, truth.

Hitherto the prestige of mechanics as a secure body of knowledge had been upheld by its presentation as a deductive science. Like Euclidean geometry, it proceeded from a few axioms (laws and principles) from which all else was deduced rationally according to mathematical logic. The axioms themselves were supposed to be truths of nature, a position that was difficult to uphold as the role of force as an axiom was challenged. By the 1870s, this uncertainty over the basis of mechanics was compounded as the status of Euclidean geometry itself became unexpectedly insecure, a development that had far-reaching epistemological implications.

Until the mid-nineteenth century the equivalence of mathematical and physical reasoning went unquestioned. Despite a move towards rigour, many mathematicians, especially in Britain, continued to accept physical arguments which endowed space with physical properties such as incompressibility, as valid means of proof. Thus space (defined as a collection of points), physical space (our perception of the space we inhabit), and geometry (a space with particular metric properties), were widely assumed to be one and the same thing [25]; it would have been anachronistic to distinguish the three concepts. Our experience of physical space guaranteed the truth of Euclidean geometry. Conversely, the logical structure of Euclidean geometry assured that generalising about physical space beyond our direct experience was a valid means of knowledge construction [26].

Thus Riemann's ideas of geometry challenged not only the belief that Euclidean was the one true geometry, but also the certainty that science could arrive at a universal physical truth that transcended experience [27,28]. They became widely known after Helmholtz's promotion of them in 1868, and William Kingdom Clifford's 1873 translation of Riemann [29,30]. Geometry could no longer be founded on self-evidently true axioms that were not open to physical questioning. Instead, geometries were based on metrical considerations that could be empirically investigated. The empiricist response was to consider the axioms of both mechanics and geometry as true only to the extent that they were supported by experience until now.

Gravitational action at a distance had been adopted by Newton as an empirically verified fact, standing in until such time as a more satisfactory and fundamental theory of locally acting forces could be developed. On astronomical scales, where the size of celestial bodies could be considered insignificant compared to the distances between them and constraints were not important, the assumption of action at a distance in Euclidean space according to a central force governed by an inverse square law worked well. Through the eighteenth century, driven by a desire to reduce the discrepancy between astronomical observations and Newtonian mechanics, mathematicians such as Euler, d'Alembert, Lagrange and Laplace had developed techniques of approximation and variation which allowed them to handle many-body problems. These efforts culminated triumphantly with the publication of Laplace's *Traité de mécanique celeste* in 1798-1825 [31].

However, as well as being an unexplained and hence occult force, the assumption of action at a distance for short range interactions forced a choice. The dominant, Laplacian, tradition considered short range – and even impact – forces as action at a distance between "molecules" (simple bodies) of matter, heat, light, or electricity [32]. It retained Newtonian mechanics by treating the molecules between which forces act as mathematical points without extension. This avoided the possibility of collisions between bodies. It permitted the development of potential theory by Laplace, Poisson, Green and Gauss. However, it seemed physically unrealistic at microscopic scale. The alternative, originally Newtonian option, was to allow the molecules a finite size but acknowledge that this created a lower distance bound to the validity of the inverse square law since the molecules were generally deemed impenetrable. Lazare Carnot and Thomas Young were among those who followed this option [33,34]. On both views, matter was primary: force was seated in individual material bodies with an ontological distinction between matter and the void over which the distance force acted.

The mathematical investigation of vortices took place in the context of the rise of field theories through the nineteenth century. This was driven by efforts to replace the unexplained action at a distance by contiguous forces acting locally and to resist materialism. It is no accident that the materialist, Laplacian, tradition had flourished in post-revolutionary France, while field theories developed in the context of the German idealist preoccupation with the relation between mind (or spirit) and matter, and the theological commitments of such British mathematicians and natural philosophers as Faraday, William Rowan Hamilton, Kelvin and Maxwell [3,22,32]. Field theories viewed the space between bodies as primary, propagating force through continuous action between neighbouring elements. Physically, the space was either comprised of a continuous distribution of local forces (as in Faraday's magnetic field) or permeated by ether whose state and motion effected the forces (as in Maxwell's "On Physical Lines of Force" [12]). Ideally, the space might also be the seat of action of spirit, free will, or God (as in Riemann's ether field or Balfour Stewart and Tait's *The Invisible Universe* [27,35]). The development and exploitation of field theories depended on adopting a deductive approach to science, recently promoted by John Herschel and William Whewell, using advanced mathematics to deduce testable observational predictions from the more abstract theory.

These mathematical field theories represented a more profound conceptual shift than just that in location of force. Field theorists turned their backs on attempts to devise explanations based on invisible molecules for phenomena such as heat and light. Instead, they claimed, the proper domain of physics was the observable and measurable. However, this did not necessarily mean that they embraced scientific naturalism or were materialist. The use of Lagrangian mechanics and partial differential equations enabled them to treat phenomena at the level of physical experience [36]. But underpinning the move was the assumption of a substratum of continuous physical space – arguably established by the Creator - whose state and processes were described by macroscopic differential equations which defined the observables. The physical nature of the substratum was irrelevant provided it gave rise to equations that agreed with observation, leaving ample scope for heuristic models and mathematical analogy for the processes in the space. They were aware of, but generally ignored, the questions Riemann had raised over the continuity of space [25,27] although they exploited his suggestion of multiply-connected space; their use of differential equations was predicated on continuity and the idea that force represented a continuous change in some state. In the 1850s Kelvin, Rankine, and other North British physicists defined that state as "energy" by drawing on potential theory and Hamilton's formulation of Lagrangian mechanics, and appropriating to their own needs the distinction Helmholtz had established between quantity and intensity of force; "quantity of force" became "energy" [3,22,37 p617-619]. For them "dynamical" connoted systems that could be analysed as matter in motion in a continuous medium using Lagrangian methods with energy as the foundational concept. By the 1860s and 70s energy was replacing force as the fundamental, conserved, basis of physics. In 1870 Clifford considered that for the majority of scientists, "Force is regarded as the great fact that lies at the bottom of all things," but that this view was being superseded [38 p121]. A few years later Tait opined that, "in another century [force] will probably have lost all but a mere antiquarian interest." [39 pxiv].

As a major proponent of the new energy programme in physics, Tait further claimed that "[energy] has been shown to have as much claim to objective reality as matter has," because it was conserved [39 p4]. But his view was still contentious. Clifford was one of those who objected. He disputed Tait's assertion, that matter and energy were objectively real because they were conserved, on philosophical grounds and pointed out that our knowledge of the conservation of both was limited by the accuracy of experiment. Claims of absolute conservation rested on an unjustified extrapolation from macroscopic experience. Similarly, theories of matter based on the ether rested on a mathematical fiction, a continuous perfect fluid, extrapolated from macroscopic imperfect fluids. Matter and energy, he concluded, were "complex mental images", of phenomena that might possibly result purely from the curvature of space – an idea that did not gain currency at the time but foreshadowed Einstein's general relativity [38 quote on p287].

Possibly the North British physicists felt able to ignore Clifford's suggestion of curved space since Helmholtz had provided them with an empirical argument for doing so. By insisting on the "axiom of free mobility", the hypothesis that all spatial figures may be translated or rotated in space without changing their form, he restricted physicists' (but not mathematicians') consideration to geometries of constant curvature. Then, as Riemann had pointed out, as far as they go, astronomical observations showed that astronomical space could not be distinguished from Euclidean space [27,29]. Certainly, in 1876, Tait referred to Helmholtz and Riemann when asserting that physical space was necessarily three-dimensional, and consigning enquiry into its nature to "purely metaphysical, and therefore of necessity absolutely barren, speculation" [39 p4]. He did concede, though, that space might have other properties in other parts of the universe.

Thus the background against which the North British physicists developed the mathematics of vortices was one in which any "force", and certainly action at a distance forces, were to be avoided if at all possible, "energy" was valid in their eyes but not yet widely accepted, and space was three dimensional and Euclidean at least locally but possibly not beyond the reaches of current observation. Space was necessarily continuous, and matter of any finite size that was different in nature from space was problematic. Helmholtz's background was slightly different, rooted in Laplacian forces between bodies and Newtonian mechanics. But for all, the impact of these debates is evident in the decisions they made in developing the mathematics of vortices, discussed in the next section.

## 3. Hydrodyanamics and Ethereal Vortex Theories

Unlike the North British who developed his work, Helmholtz's 1858 pioneering investigation of vortices in perfect fluids originated in a concern with real fluids - in acoustical questions - when he investigated the use of Green's theorem to analyse the forms of motion induced in fluids by friction [1]. As in his earlier "On the Conservation of Force" [20] his arguments were based on action at a distance forces and potential theory. His belief in bodies as fundamental showed as he traced the motion of "particles" of fluid using material coordinates that travelled with them. He conceived differences in potential physically as the gain or loss of *vis viva* (we would say the work done) in moving from one position to another. Following Kirchhoff, and like Stokes, of whose work he was initially unaware, Helmholtz decomposed the infinitesimal motion of a volume element of fluid into a translation, three orthogonal dilations, and a rotation. He devised the basic definitions of vortex line and vortex filament (the tube formed by the vortex lines through every point around an infinitesimally small closed curve). Then, making widely used, but physically questionable, assumptions that the surrounding fluid acts as a normal pressure on the surface of a fluid element, that the density was essentially constant, and that Newton's laws could be applied unproblematically to the assemblage of points comprising the "particle", he showed that each vortex line or filament always comprised the same elements of fluid. Fluid that was not rotating initially would not start, and fluid that was rotating to start with would never stop; it stayed inside the vortex filament and moved with the filament. As Lamb subsequently emphasised, this was a conclusion about particles of fluid, not about portions of space [40]. Thus, each filament had a unique stable identity and could be treated as a body. Moreover, vortex lines would either form closed rings or terminate on the boundary of the fluid. Either way they were permanent, and their "strength", the product of their cross-section and their angular momentum, was invariant: the strength was the same everywhere around a vortex ring and remained constant as the ring moved through the fluid.

Throughout his paper, Helmholtz used an analogy between hydro- and electrodynamics, "to give a briefer and more vivid representation" [2 p487,41]. He pointed out that the lines of fluid motion were arranged in the same way as the lines of magnetic force around an electric current flowing along the axis of the vortex, and worked from his previous understanding of electromagnetism to develop physical arguments and gain an intuitive grasp of vortex behaviour [1]. He first utilised the analogy when attempting to find the velocity distribution within a mass of fluid, if the distribution of vorticity was known - the resultant equations had the same form as those for the magnetic field due to a current-carrying wire. In the following sections he exploited the analogy, picking up on Riemann's ideas of multiply-connected space and specifying functions by their boundary conditions, published the previous year [42]. Reasoning from the observed properties of electricity and magnetism Helmholtz investigated the fluid flow in the multiply-connected space around the vortices and the action of vortices on each other. He suggested that the velocity potential outside the vortices must be a many-valued harmonic function, and gave $\tan^{-1}(x/y)$ as an example of such a function. His explanation is not entirely clear, but appears to be that, unlike the gravitational case, the work done in moving between two points in the space (mass of fluid or magnetic field) *does* depend on the path taken. Holes (vortex filament or current-carrying wire) represent singularities in the space and the potential increases each time the path circles a hole. In such cases Green's theorem did not apply. The velocity in the fluid could not be specified by the boundary conditions. Except in a few simple cases it could not be found analytically, but could, by analogy, be represented by the magnetic force around current-carrying wires. Although Helmholtz did not comment on this, the analogy worked only if the vortices extended from one boundary to another or formed closed loops, just as current-carrying wires either linked conductors at different potentials or formed closed loops. This simple analogy to wires helps to explain his conclusion that if the vortices did not end on boundaries, they must be closed. It parallels Faraday's distinction between matter, which could initiate lines of magnetic action, and space, which could not [43]. The conclusion is unfounded topologically and has puzzled Epple [44], although it was widely accepted as valid at the time.

To exploit the analogy further, Helmholtz explicitly assumed that vortices that had both ends on the same surface could be joined up into rings in imagination outside the surface, the space considered being enlarged to encompass the entire ring. To formulate tractable simple cases, Helmholtz imposed the condition that "the rotation of the fluid elements takes place only in known surfaces or lines," whose form remained unchanged during the motion. Under these conditions, he analysed a vortex sheet (discussed below), and the action of two vortices on each other: if they were parallel filaments he concluded that they circled each other around a common "centre of mass"; if they were two rings close to each other with parallel vorticity and with their centres on the same axis perpendicular to the rings, they leap-frogged one another, the front one slowing and enlarging while the back one shrank and speeded up until it passed through the centre of the first and, in its turn, slowed and enlarged [2]. Antiparallel vortex rings approached each other, ever enlarging and slowing asymptotically.

Maxwell was apparently unaware of Helmholtz's vortex paper and electromagnetic analogy when he embarked on his own electromagnetic theory by using vortex lines as an analogy for the lines of magnetic force in "On Physical Lines of Force", published in instalments in 1861-62 [12]. Unlike Helmholtz, who likened vortex lines to current-carrying wires and worked from his previous understanding of electromagnetism, Maxwell did the reverse [4,45]. He used his knowledge of stresses and strains in continuous fluids to generate a mathematical theory of electromagnetism. Using mainly geometric and physical arguments, he devised a fluid-based model of the magnetic field as an array of parallel vortices. Late in the paper, he justified his choice by asserting that the rotation of the plane of polarised light by a magnetic field showed that "the cause of the magnetic action on light must be a real rotation going on in the magnetic field" [12 p88], a suggestion that was due to Kelvin. His next step, guided by imagining real, solid, mechanical systems, was to annul the discontinuity in velocity between contiguous vortices by inserting imaginary "idle wheels" between them. The displacement of the idle wheels came to represent electric current [12]. This step exemplifies a further contrast with Helmholtz's work, arising from Maxwell's different purposes for the analogy. Where Helmholtz used a "sheet" of parallel vortices specifically to show the possibility of tangential velocity discontinuities in the fluid on either side of the "indefinitely" thin sheet, and was not bothered by the issue of discontinuity in velocity between contiguous vortices within the sheet, Maxwell insisted on continuity. Whereas Helmholtz's aim was to understand real fluids ultimately composed of particles acting at a distance, Maxwell was interested in vortex-like processes in space; both his beliefs about space and known observations of electricity and magnetism suggested continuity.

Maxwell emphasised the role of his model at the outset of the paper: "I propose now to examine magnetic phenomena from a mechanical point of view, and to determine what tensions in, or motions of, a medium are capable of producing the mechanical phenomena observed. If, by the same hypothesis, we can connect the phenomena of magnetic attraction with electromagnetic phenomena and with those of induced currents, we shall have found a theory which, if not true, can only be proved to be erroneous by experiments which will greatly enlarge our knowledge of this part of physics" [12 p162]. Maxwell's usage of "hypothesis" here reflects the fact that he was doing something methodologically new and going beyond his contemporaries' understanding of "model" [46]. The hypothesis could suggest experiments for testing. Provided the theory correctly described known observations, its hypothetical basis could not be challenged except by extending the realm of observation. Being purely a heuristic model, Maxwell did not concern himself with the questions of permanence and stability of vortices that later teased Kelvin and Helmholtz, and that his consideration of strain patterns might have exacerbated. Nor was he interested at this point in developing the mathematics of vortices as such; the paper came to no new conclusions about vortices. He subsequently discarded the model, which had served its purpose in leading him to a mathematical theory of the electromagnetic field [15]. Nevertheless, his use of the analogy had extended the domain of application of vortex mathematics to a vitally important new area of physics and boosted the confidence of British physicists in the power of hydrodynamical models.

Kelvin was one such physicist, with a long history of using hydrodynamical, though not vortex, models [47]. At the same time that physicists were rejecting molecular explanations of light and electricity, as noted in the previous section, they were re-introducing them as explanation for the macroscopic properties of gases and the spectra of the elements. Such phenomena required more complex molecules or atoms than the simple hard bodies of Laplacian physics. Kelvin turned to the ether for his molecular model. Initially, he suggested privately to Stokes that discrete atoms might originate in some form of inhomogeneity in the continuum [48 p330]. But in 1867 Tait's experiments with smoke rings drew his attention to Helmholtz's hydrodynamics which, Kelvin realised, had laid the groundwork – permanence and individual identity - for an atomic model based on ethereal vortex rings [17,49].

Helmholtz never suggested a vortex atom theory himself and continued to work with central action at a distance forces, remaining aloof from Kelvin's efforts to use continuum mathematics, energy, and impact forces. Kelvin adapted Helmholtz's work to his own interests, as he had done in the case of potential noted above. In 1867 he emphasised in an excited letter to Helmholtz, that a vortex atom in a perfect fluid ether would not only, "be as permanent as the solid hard atoms assumed by Lucretius", but permanence also allowed different types of structured atoms, for "if two vortex-rings were once created in a perfect fluid, passing through one another like links of a chain, they never could come into collision, or break one another, they would form an indestructible atom; every variety of combinations might exist." [19 vol1 p514-5]. Vortex rings could associate with each other in pairs or groups and had definite modes of vibration which might explain the spectra of the elements [49]. Kelvin followed up his initial suggestion with a long mathematical analysis of both rotational and irrotational motion in a perfect fluid, using energy and impulse as his basic concepts [50]. His aim was "to illustrate the hypothesis, that space is continuously occupied by an incompressible frictionless liquid acted on by no force, and that material phenomena of every kind depend solely on motions created in this liquid" [50 p217]. Thus, he was more explicit than Helmholtz about the distinction between "space" (irrotational motion) and "bodies" (rotational motion) and devoted much attention to the movement of "bodies" through irrotational fluid. His main hydrodynamical insight was that Helmholtz's theorems could all be encompassed within a single theorem: that the circulation of the vortex remains the same throughout motion, where he defined circulation as the line integral of the tangential velocity component around a closed curve in the fluid, equivalent to the flux of vorticity through a surface defined by the curve [17,50].

The enthusiasm of many British mathematical physicists for Kelvin's suggestion originated in its alleged freedom from assumptions, its reliance on a sophisticated mathematical programme with which they were already asserting their authority in other areas of physics, and their religious convictions. Maxwell promoted the idea strongly in his address to the Mathematical and Physical Section of the British Association in 1870 and in his influential *Encyclopaedia Britannica* article, "Atom", in 1875. The theory appealed to him because, he alleged, unlike other atom theories, it was closed and free from arbitrary assumptions. It rested on "nothing but matter and motion" and included "no central forces or occult properties of any kind," [51 p223]. Thomson similarly emphasised that his investigation of vortex rings was "almost entirely kinematical" and did not depend on assumptions about forces [8 p2]. Responding to the debate over force, they both felt that, as Truesdell later remarked, "a kinematical result is a result valid forever, no matter how time and fashion may change the "laws" of physics," [52 p2]. The complexity of the mathematics was an attraction for Maxwell's student, Donald MacAlister: "the work of [mathematical] deduction is so difficult and intricate that it will be long before the resources of the theory are exhausted" [53 p279]. Maxwell concurred: "the difficulties of this method are enormous, but the glory of surmounting them would be unique." He asserted that the ether, "has no other properties than inertia, invariable density, and perfect mobility, and the method by which the motion of this fluid is to be traced is pure mathematical analysis," [54 p472]. However, despite Maxwell's disclaimer, the origin of vortex atoms *was* mysterious. As Thomson noted, "This theory cannot be said to explain what matter is, since it postulates the existence of a fluid possessing inertia" [8 p1]. Nor did it explain the origin of the vortices. In *"Atom"* Maxwell glossed over the issue: "When the vortex atom is once set in motion, all its properties are absolutely fixed," he stated, without asking how they were set in motion [54 p471]. However, he was more explicit later in the article. Reviving John Herschel's observation that atoms were like "manufactured articles", he argued that the uniformity of the properties of atoms through time and space (whatever the underlying atomic model) was beyond the power of science to explain: "We have something which has existed either from eternity or at least from times anterior to the existing order of nature…. Science is incompetent to reason upon the creation of matter itself out of nothing," [54 p482]. Maxwell was writing just after Tyndall's "Belfast Address" to the British Association for the Advancement of Science, which strongly promoted materialism. In opposition, Maxwell pointed to the limitations of science, and deliberately left open the possibility of a manufacturer or Creator. The following year Tait and Stewart exploited this possibility of an intelligent design argument for God in *The Unseen Universe* [35], basing their account of matter on vortex rings. Once again Clifford was an energetic opponent. He pointed out that the theory supposed that different ring topologies represented different elements. So, as well as creating an ether, divine intervention had to have set in motion (at least) 63 different forms of vortex rings; the theory was not as economical as claimed [38 p296]. Vortex atom enthusiasts ignored his views.

Thus, the observation that atomic properties were invariable underpinned the attraction of vortex atoms to many British physicists both theologically and mathematically and guided their investigations. Although he proved the permanence of single vortex rings, Helmholtz had considered only vortices whose form did not change through time. On Kelvin and Tait's supposition that difference in topological form explained the physical properties of different types of molecule, it was crucial to establish the stability of form of rings within a space that contained multiple rings, that might interact with each other as well as with the surrounding irrotational ether. As early as his 1869 paper, Kelvin addressed the question of the stability of the flow of fluid outside vortex filaments. By considering the complex, but topologically invariant, vortex rings as "bodies" he treated the space between them as multiply-continuous and tackled the problem of modifying Green's theorem so it could still be used to analyse the flow. His ongoing use of the term "multiply continuous" in this and later papers is significant. He was well aware that Tait had translated Helmholtz's *zusammenhängend* as "connected" and of Tait's reasons for that decision. But Kelvin's focus was on flow in the space between the holes, and here continuity was essential to his mathematical approach. In debates with Tait and Helmholtz he insisted that the same held physically for a perfect fluid [44], a claim that Clifford challenged as observationally unjustified, as discussed in the previous section.

Kelvin's involvement with nascent ideas of topology is evident as he struggled for definitions: "I shall call a finite portion of space n-ply continuous when its bounding surface is such that there are n irreconcilable paths between any two points on it," he declared [50 p44], using a path definition akin to later ideas of homotopy. However, he went on to work, instead, with the physical idea of "stopping barriers", virtual membranes that prevented the motion of the fluid, reducing the space to simply continuous and removing the ambiguity in the value of the velocity potential. This intuitive physical notion, related to what we know as homology, fitted well with Riemann's or Helmholtz's ideas of cutting surfaces [44]. In the end, Kelvin conceded in a footnote that Helmholtz's surfaces were more physically useful than his own path definition [50 p55]. Assuming sufficiently smooth functions, Kelvin showed how to modify Green's theorem and determine the flow in the irrotational fluid by including terms representing integration along a path from one side of each barrier to the other. These path integrals became the "circulation" with which he re-cast Helmholtz's theorems. However, his subsequent attempts to use energy considerations to show that systems of vortex rings were stable under perturbation met with little success [4,44]. Even ten years later, the only case for which he had mathematical proof of stability was the very simplest – a single columnar vortex – for which he showed that a small periodic deformation of the surface propagated around and along the column with constant amplitude [55]. He claimed, using physical reasoning rather than mathematical proof, that the motion of a perfect fluid in a finite container would be both steady and stable if the kinetic energy of the fluid was a maximum or minimum [56]. However, the kinetic energy of a vortex ring was, instead, a "minimax" or saddle point, so energy considerations were unable to determine the stability of vortex atoms. Kelvin's failure highlights his commitment to a perfect fluid analogy, and its limitations: Moffat [57] suggests that had Kelvin employed a magnetic analogy in this case, pre-empting the development of magnetohydrodynamics and using the magnetic induction equation (which was known to him) in place of the vorticity equation, the outcome might have been different, since such stable, knotted magnetic structures do exist.

## 4. Debates over continuity

Lack of demonstrable stability was not the only problem that came to beset the vortex atom theory. Many issues originated in concepts of continuity derived from extrapolation of macroscopic physical experience to the microscopic, or in transitioning from purely kinematical mathematical analysis of points and lines to dynamical analysis of physical bodies endowed with mass and energy. As Truesdell later observed, while a purely kinematic theory might be true, it is "inadequate to represent more than a highly restricted and indeed dynamically degenerate class of phenomena," of little relevance to physics [52 p72]. However, application of Newton's laws involved an unwarranted extension to infinitesimal elements of mass if the fluid was a true continuum, and authors vacillated inconsistently between treating elements as "particles" or "points". Osborne Reynolds criticised Lamb: "Having rightly defined fluids as being such 'that the properties of the smallest portions into which we can conceive them divided are the same as those of the substance in bulk,' he proceeds to reason about a particle as though it were a discrete quantity, the position of which is defined by some point," and this had led to errors in Lamb's argument that particles that lie in a bounding surface at one time always lie in it [58 p343].

Kelvin had hoped that ''A full mathematical investigation of the mutual action between two vortex rings . . . will become the foundation of the proposed new kinetic theory of gases'' [49 p16]. In this he was disappointed. Vortex rings might provide a plausible explanation for molecular spectra, but they came into increasing conflict with the kinetic theory. This theory regarded molecules as hard elastic balls whose motion explained the thermal properties of the gas. It was already subject to criticism, for the equipartition theorem, which stated that in thermal equilibrium the energy of the molecules should be shared equally between its various forms, predicted a ratio of specific heats for a diatomic molecule of 1.33, whereas the observed value was 1.44. The situation for vortex molecules was even worse, for continuity implied an infinite number of degrees of freedom. Maxwell pointed out this problem in 1877,

> *It will not do to take a body formed of continuous matter endowed with elastic properties, and to increase the coefficients of elasticity without limit till the body becomes practically rigid. For such a body… has an infinite number of degrees of freedom…. Such atoms would soon convert all their energy of agitation into internal energy, and the specific heat of a substance composed of them would be infinite* [17 p58].

A similar problem underlay Kelvin and Tait's analysis of the movement of bodies (vortex atoms) through the surrounding irrotational fluid [59], an analysis essential to a vortex theory of gases. Using Lagrange's equations of motion, they characterised the motion in terms of kinetic and potential energies and generalised coordinates, avoiding any mention of forces. Hicks praises this step as, "the next great advance in theory" after Helmholtz's 1858 paper [7 p60]. The method was favoured by the North British physicists for analysing situations where they wished to postulate an underlying mechanism. Maxwell pointed to the difference from astronomy where Newtonian theory worked well, for its role was purely to account for the observations. However,

> *when we … believe that the phenomena which fall under our observation form but a very small part of what is really going on in the system, the question is… what is the*

> *most general specification … consistent with the condition that the motions of those parts of the system which we can observe are what we find them to be? It is to Lagrange … that we owe the method which enables us to answer this question without asserting either more or less than all that can be legitimately deduced from the observed facts"* [60 pp781-2].

But Kelvin and Tait's application of the method to vortices was widely attacked, for example by Ludwig Boltzmann, John Purser (Professor of Mathematics at the Queen's College, Belfast), and Carl Neumann [61]. The problem was the lack of any proof of the validity of the equations in situations where there were infinite degrees of freedom, which were not recognised in the generalised coordinates.

In the second edition of their *Treatise* Kelvin and Tait re-affirmed their commitment to the method by promoting a process they termed "ignoration of coordinates" [59 p327]. It was essentially the same as that developed two years earlier by the Cambridge mathematician Edward Routh in his Adams Prize essay on stability of motion, leading to a "modified Lagrangian"[62]. The ignorable coordinates (later termed cyclic) were those that did not appear explicitly in the kinetic or potential energies of the system, and hence not in its Lagrangian. As realised by Louiville, such ignorable coordinates imply that their conjugate momenta are constant [23 p209]; Routh's modified Lagrangian included the constants explicitly, but they were then eliminated in the equations of motion, leaving just the explicitly appearing coordinates; the non-explicit coordinates were "ignored". Hicks, in his report, accepted Kelvin and Tait's solution to the problem and credits them with, "transforming the methods of hydrodynamics" [7 p60]. But Neumann pointed out that this move still depended on the ability to define the coordinates. Kelvin and Tait had not considered "the transition from a finite to an infinite number of degrees of freedom and the difficulty of identifying the coordinates to be ignored" [61 p107 author's translation]. Kirchhoff had derived the equations directly from Hamilton's principle, and Neumann deemed this more acceptable.

Ultimately, lack of assured stability was one reason the vortex atom programme foundered. However, Kelvin's debates with Stokes about stability reveal the importance of differing assumptions about continuity to the mathematical development. As early as 1859 Stokes had cautioned Kelvin that the motion of a perfect fluid around a solid body was unstable [4]. He suggested that surfaces of discontinuity would form in the wake of the body, implying a loss of *vis viva* and resistance to the motion. Similarly, a surface of discontinuity would be formed, he believed, in a fluid passing an edge from one container to another. Kelvin disagreed; Lagrange's theorem, he argued, showed that the motion was irrotational with no eddies. Stokes responded that Lagrange's theorem assumed continuity as a necessary condition, so was inapplicable in the presence of a surface of discontinuity. Twenty years later Stokes reflected that his own ideas of continuity were an extrapolation from the behaviour of fluids with low viscosity. Similarly, he had assumed a perfectly sharp edge. Beginning with a real fluid and a real edge, the outcome for perfect fluid and edge depended on the order in which one took the limits as the viscosity and the radius of curvature tended to zero. Moving to a zero-viscosity perfect fluid first resulted in continuous irrotational flow regardless of the radius of curvature. But moving to a perfect edge first resulted in a trail of vorticity from the edge that becomes a vortex sheet in the zero-viscosity limit.

Helmholtz returned to his investigation of discontinuity surfaces in 1868 [63], representing them, as ten years earlier, by vortex sheets. He argued that the sheet must move parallel to itself with a velocity the mean of those on either side of the sheet. Physical intuition suggested that in a real fluid with viscosity, the fluid on either side of the sheet would be gradually set in motion and the sheet would grow to a set of finite-sized vortices. Most significantly, experiments showed that vortex sheets were highly unstable; Helmholtz claimed that any perturbation would grow, leading to "a progressive spiral unrolling of the corresponding portion of the surface which portion, moreover, slides along the jet" [4 p161].

The interest in vortices of Stokes and Helmholtz, and consequently their approach, differed markedly from that of Kelvin and Maxwell. Stokes and Helmholtz were interested in the behaviour of real fluids. They studied perfect fluids as an intermediate step on the way to real fluid behaviour. Kelvin and Maxwell, on the other hand, were interested in energy-based dynamics of a perfect fluid ether as a model that might or might not be real. A presupposition of continuity of motion was a necessary condition for their Lagrangian methods.

## 5. Conclusion

By the time of Maxwell's death in 1879, the basic laws of vortices in a perfect fluid in three-dimensional Euclidean space had been established, as had their importance to physics. For one, dominant, group of British mathematical physicists they were supremely important. Epple has pointed out that this group was pursuing what Scholz has called "heteronomous mathematics," [44,64]. Their motivations, and many of their arguments, were physical, but they were developing new mathematics rather than applying existing methods. However, tensions were apparent throughout, with no simple directionality. Sometimes assumptions required for mathematical tractability conflicted with physical theory as in the assumption of a perfect fluid whose action on an element of the fluid could be represented by a normal pressure, which James Jeans characterized as, "an abstract ideal which is logically inconsistent with the molecular constitution of matter postulated by the Kinetic Theory," [65 p29]. At other times the same assumptions were made for physical or theological objectives, such as ensuring the permanence of atoms, but ran into mathematical problems, as in Kelvin and Tait's attempts to apply Lagrange's equations of motion to vortices.

Osborne Reynolds was clear that the origin of the tensions lay in the impossible attempt to match abstract theory with reality. The fluid dynamicist,
> *has to begin by apologising for his fundamental assumptions as being obviously contrary to facts, and after carrying his readers through most difficult and complex mathematics, he has again to apologise for his conclusions, which are in general contrary to experience… The only real success, has been obtained by the rigorous development of the theory of the motions in a perfect fluid… regardless of whether or not these motions take place in actual fluids."* [58 p343].

Vortex atom theory in Britain was pursued by an elite group of Cambridge-trained mathematicians and physicists. It provided a theological and physical context for exploiting their own mathematical training, and a means of exerting authority over less mathematically-trained physicists. But ultimately, they were unable to resolve what MacAlister deemed the, "true paradox that one must be something of a mathematician to understand what a *real* theory the Vortex-theory is" [53 p279 MacAlister's emphasis].

## Acknowledgements

I am grateful for the perceptive comments of two anonymous referees, which have enhanced this paper.